\documentstyle[12pt]{article}

\newcommand{\supp}{\mathop{\rm supp}\limits}

\begin{document}


\begin{center}

{\bf COMPACT SETS IN BIDE - SIDE GRAND LEBESGUE SPACES, \\

\vspace{3mm}

 WITH APPLICATIONS.}\\

\vspace{3mm}

{\sc Eugene Ostrovsky, Leonid Sirota}\\

\vspace{3mm}

 {\it Department of Mathematic, HADAS company, \\
56209, Rosh Ha Ayn, Hamelecha street, 22; ISRAEL;} \\
 E - mail: galo@list.ru, \ eugeny@soniclinx.com \\

\vspace{3mm}

{\it Bar \ - \ Ilan University,  59200, Ramat Gan, ISRAEL;} \\
 e - mail: sirota@zahav.net.il \\

\vspace{3mm}

              {\sc Abstract.}\\
\end{center}

\textwidth = 12 cm

\begin{verse}

\normalsize
\hspace{5mm} In this article we find some sufficient \\
\hspace{5mm} conditions for the set in the Bilateral \\
\hspace{5mm} Grand Lebesgue Space to be compact set.\\
\hspace{5mm} We consider applications into numerical \\
\hspace{5mm} methods and  in the basis problem .\\

\textwidth = 15 cm

 \vspace{3mm}


\vspace{3mm}

Secondary (2000) 34A34, 65M20, 42B25, 35A25 . \\

\end{verse}

\vspace{3mm}

{\bf 1. \ Introduction. Statement of problem.}\\

\vspace{3mm}

 Let $ (X, \Sigma,\mu) $ be a measurable space with
non \ - \ trivial measure $ \mu: \ \exists A \in \Sigma, \mu(A) \in
(0,\mu(X)). $ We will assume that either $ \mu(X) = 1 \ $ or $
\mu(X) = \infty $ and that the measure $ \mu \  $ is $ \sigma - $
finite and diffuse on the sets of positive $ \ \mu \ $ measure:
 $ \forall A \in \Sigma, 0 < \mu(A) < \infty \
\exists B \subset A, \mu(B) = \mu(A)/2. $  Define as usually for
arbitrary measurable function $ f: X \to R^1 \ $
$$
|f|_p = \left(\int_X |f(x)|^p \ \mu(dx) \right)^{1/p}, \ p \ge 1;
$$
$ L_p = L(p) = L(p; X,\mu) = \{f, |f|_p < \infty \}. $
 Let $ a = const \ge 1, b = const \in (a,\infty], $ and let $ \psi = \psi(p) $
be some strictly  positive on the {\it closed } interval $ [a,b] $, including
the cases

$$
 \psi(a) \stackrel{def}{=} \lim_{p \to a + 0} \psi(p) = \infty
$$
 and

$$
\psi(b) \stackrel{def}{=} \lim_{p \to b \ - \ 0} \psi(p) = \infty,
$$
continuous on the {\it open} interval $ (a,b) $ function. \par

 The set of all those functions we will denote $ \Psi: \ \Psi = \Psi(a,b) =
\{ \psi(\cdot) \}. $

We define in the case $ b = \infty \ \psi(b-0) = \lim_{p \to \infty}
\psi(p). $ \par

\vspace{3mm}

{\bf Definition 1.} \\

\vspace{3mm}

 Let $ \psi(\cdot) \in
 \Psi(a,b). $  The space $ BSGL(\psi) = G(\psi) = G(X,\psi) =
G(X,\psi, \mu) = G(X,\psi,\mu, a,b) $ (Bide \ - \ Side Grand Lebesgue
Space) consist on all the measurable functions $ f: X \to R $ with
finite norm
$$
||f||G(\psi) \stackrel{def}{=} \sup_{p \in (a,b)} \left[|f|_p/\psi(p) \right].
$$

The introduced spaces are some generalization of the
so-called {\it Grand Lebesgue Spaces} (see
 \cite{Fiorenza1}, \cite{Fiorenza2},  \cite{Iwaniec1},  \cite{Iwaniec2},  \cite{Jawerth1} etc). For example, the space $ L^{b)}, b \in [1, \infty) $
(in the notations  \cite{Fiorenza1}, \cite{Fiorenza2},  \cite{Iwaniec1})
 coincides with our space $ G(\psi_b(p)), $ where by definition

$$
\psi_b(p) = (b \ - \ p)^{\ - \ 1/b}, \ p \in (1, b).
$$

 They are also rearrangement invariant (r.i.) spaces, and was studied as r.i.
spaces  (find its fundamental function,  calculated its adjoin  spaces, obtained 
different imbedding and convolution theorems etc.) in the papers, e.g.,   
\cite{Ostrovsky3},  \cite{Ostrovsky5}. Moreover, the $ G(\psi) $ spaces are the so \ - \ called moment rearrangement invariant (m.r.i.) spaces. They satisfy the Fatou property etc. \cite{Ostrovsky7}.

{\bf We investigate in this article some  conditions for arbitrary subsets
of $ G(\psi) $ spaces to be compact set in this space or in some its closed
subspace. We intent also to consider some applications: in the theory of numerical methods and investigate existence of basis in these spaces.} \par

\vspace{2mm}

 Note that the $ G(\psi) $ spaces are the particular case of interpolation spaces (so-called  $ \Sigma - $ spaces) (see \cite{Astashkin1}, \cite{Astashkin2}, \cite{Capone1},\cite{Carro1}, \cite{Jawerth1}).\par
 But we hope that the our direct representation of these spaces (definition 1) is more convenient for investigation and application.\par
 In the case  $ \mu(X) = 1 \ G(\psi) \ $ spaces appear in the article
\cite{Kozatchenko1}. See for  detail \cite{Ostrovsky2}, chapter 1.   \par

  These spaces are used, for example, in the theory of probability
 (\cite{Kozatchenko1}, \cite{Ledoux1}, \cite{Ostrovsky1}, \cite{Ostrovsky2}, \cite{Ostrovsky3}, \cite{Ostrovsky4}, \cite{Ostrovsky5},\cite{Ostrovsky6},
\cite{Talenti1} etc.), theory of PDE (\cite{Fiorenza2},  \cite{Iwaniec2}),
 functional analysis \cite{Davis1}, \cite{Fiorenza1},  \cite{Jawerth1},
theory of Fourier series (\cite{Ostrovsky5}), theory of martingales
\cite{Ostrovsky1}, \cite{Ostrovsky4}, theory of approximation
 \cite{Ostrovsky7} etc.\par

 The article is organized as follows. In the next section we reproduce
some affirmations from the general theory of BSGL spaces. The third section
contains the main results. In the last section we consider some application
of obtained in section 3 results.\par

\vspace{4mm}

{\bf 2. Auxiliary facts. Subspaces. Convergence.} \\

\vspace{3mm}

 The next facts about $ G(\psi) $ spaces are proved in \cite{Ostrovsky5}.
Some assertions on the $ G(\psi) $ spaces may be obtained from the general
theory of rearrangement invariant (r.i.) spaces, see \cite{Bennet1}, \cite{Krein1},
as long as the BSGL spaces are r.i. spaces. \par

 We will say as usually \cite{Bennet1}, pp. 14 \ - \ 16 \
that the function $ f \in G(\psi), \ \psi \in \Psi $ has
{\it absolutely continuous norm and write }$ f \in GA(\psi), $
{\it if}
$$
 \lim_{\delta \to 0+} \sup_{A: \mu(A) \le \delta} ||f \ I_A||G(\psi)=0;
$$

$$
 I_A = I_A(x) = 1, \ x \in A; \ I_A = I_A(x) = 0, \ x \notin A.
$$

 A {\it family } $ F $ of elements of $ G(\psi) $ space is said to be
equi \ -  \ absolute continuous, write: $ F \in EGA, $ iff

$$
 \lim_{\delta \to 0+} \sup_{A: \mu(A) \le \delta} \sup_{f \in F}
||f \ I_A||G(\psi)=0.
$$

 We denote by $  G^o = G^o_X(\psi), \ \psi \in \Psi $ the closed subspace of
 $ G(\psi), $ consisting on all the functions $ f, $ satisfying the following condition:
$$
\lim_{p \to a+0} |f|_p / \psi(p) = \lim_{p \to b-0} |f|_p /\psi(p) = 0,
$$
in the case $ \psi(a+0)=\infty, \ \psi(b-0) = \infty; $
$$
\lim_{p \to b-0} |f|_p/\psi(p) = 0
$$
in the case $ \psi(a+0) < \infty, \ \psi(b-0) = \infty; $
$$
\lim_{p \to a+0} |f|_p/\psi(p) = 0
$$
in the case $ \psi(a+0) = \infty, \ \psi(b-0) < \infty; $ briefly:

$$
\lim_{\psi(p) \to \infty} |f|_p/\psi(p) = 0.
$$

 {\it In the case $ \psi(a) < \infty, \psi(b) < \infty $ the space $ G(\psi) $
coincides  up to norm equivalence with the direct sum $ L_a + L_b; $
which is the known Orlicz space satisfying a $ \Delta_2 $ condition. This case is not interest for us and must be excluded. } \par

 We  denote also by $ GB =  GB(\psi) $ the closed span in the norm $ G(\psi) $
the set of all bounded: $ \sup_{x \in X} |f(x)| < \infty $
 measurable functions with finite support: $ \mu(\supp \ |f|) < \infty. $  \par

The subspaces $ GA(\psi), GB(\psi), G^0(\psi) $ are closed
subspaces of the space $ G(\psi) $ and moreover coincide:

$$
GA(\psi)= GB(\psi) = G^0(\psi),
$$
 see \cite{Ostrovsky5}.\par

 Note that in the considered case $ \sup_{p \in (a,b)} \psi(p) = \infty $
the subspaces $ GA(\psi), GB(\psi), G^0(\psi) $ are {\it strong } subspaces of
the space $ G(\psi): $

$$
G(\psi) \supset GA(\psi)= GB(\psi) = G^0(\psi) \ne G(\psi).
$$
 If the metric space $ (\Sigma, \rho), $ where

$$
\rho(A_1, A_2) = \mu(A_1 \setminus A_2) + \mu(A_2 \setminus A_1), \ A_1, A_2
\in \Sigma
$$
is separable, then the spaces $ GA(\psi), GB(\psi), G^0(\psi) $ are also separable. In contradiction, the "generic" space $ G(\psi) $ is not separable.\par

 Further, let $ \psi(\cdot), \ \nu(\cdot) $ be two functions from the set
$ G(\psi;a,b). $ We will write $ \nu << \psi, $ or equally $ \psi >>\nu, $
iff

$$
\lim_{\psi(p) \to \infty} \nu(p)/\psi(p) = 0.
$$

\vspace{3mm}

{\bf 3. Main results.} \\

\vspace{3mm}

{\bf A.} In this subsection we consider arbitrary {\it bounded closed}
 subset $ S $ of the space $ G(\psi) = G(\psi; a,b): \ S \subset G(\psi). $ \par

 Since all the considered spaces: $ G(\psi), GA(\psi) $ etc. are complete metric spaces, we can and will restrict itself only by the notion of {\it sequentially compactness}.\par

\vspace{3mm}

{\bf Theorem 1.} {\it Let $ \psi, \ \nu $ be two functions from the set
$ G(\nu; a,b) $ such that $ \psi(\cdot) >> \nu(\cdot). $  Assume that: } \\

{\bf 1.}

$$
\sup_{f \in S} ||f||G(\nu) < \infty.
$$

{\bf 2.} {\it For all the values $ p \in (a,b) $ the set $ S $ is compact 
set in the space} $ L_p = L_p(X; \mu). $ \\

{\it Then the set $ S $ is compact set in the space} $ G(\psi). $ \par

\vspace{3mm}

{\bf Proof.} Let $ \{ f(n) = f(n,x) \} $ be arbitrary sequence functions
belonging the set $ S. $  Without loss of generality we conclude by virtue of condition
{\bf 1 } that

$$
\sup_n |f(n, \cdot)|_p \le \nu(p).
$$

 Let also $ \{ p(i), \ i = 1,2,\ldots \} $ be arbitrary dense in
the ordinary distance $ d(p,q) = |p \ - \ q| $ sequence of numbers in the interval
$ (a,b). $ There exists a sequence $ n(j) $ for which there exists a limit in $ L_{p(1)} $
sense:

$$
 | f(n(j), \cdot) \ - \ h(1,\cdot)|_{p(1)} \to 0, \ j \to \infty.
$$

Choosing an appropriate subsequence $ n(j(l)), $ we can see

$$
 | f(n(j(l)), \cdot) \ - \ h(2,\cdot)|_{p(2)} \to 0, \ l \to \infty.
$$

Analogously

$$
 | f(n(j(l(q)), \cdot) \ - \ h(3,\cdot)|_{p(3)} \to 0, \ q \to \infty
$$
etc. It is evident that

$$
\mu \{ \cup_{i=2}^{\infty} \{ x: h(i,x) \ne h(1,x)   \}  \} = 0,
$$
therefore we can assume that $ \forall i \ge 2  \ \Rightarrow h(i,x) = h(1,x) $
almost everywhere. \par

 Extracting the "diagonal" subsequence, indeed: $ g(1,x) = f(1,x), \ g(2,x)
 = f(n(2),x), \ g(3,x) = f(n(j(3)),x), \ g(4,x) = f(n(j(l(4))), x) $ etc., we
 conclude that for all values $ p(i) \in (a,b) $

 $$
 \lim_{m \to \infty} |g(m,\cdot) \ - \ h(1,\cdot)|_{p(i)} = 0.
 $$

 Since the sequence $ \{ p(i) \} $ is dense in the interval $ (a,b), $
we find on the basis of H\"older inequality that for all values $ p,
\ p \in (a,b) $

 $$
 \lim_{m \to \infty} |g(m,\cdot) \ - \ h(1,\cdot)|_{p} = 0.
 $$

Tacking into account the inequality

$$
|h(q,\cdot)|_p \le \sup_n |f_n|_p \le \nu(p),  \ \forall q = 1,2,\ldots,
$$
we see that $ h(1,\cdot) \in G(\nu;a,b). $ \par

 The convergence

 $$
 g(m,\cdot) \stackrel{G(\psi)}{\to} h(1,\cdot)
 $$
 it follows from a results of a paper \cite{Ostrovsky5}. \par

\vspace{3mm}

{\bf Remark 1.} It is evident that the first condition of the Theorem 1 is
necessary and that the second is not. \par
 Let us consider the following example. Choosing the function $ f = f(x) $
from the set $ G(\psi) \setminus G^o(\psi): $

$$
f(\cdot) \in G(\psi) \setminus G^o(\psi), \
$$
and such that $ ||f||G(\psi) = 1; $ we may consider the sequence

$$
f(n,x) = f(x) \cdot (1 \ - \ 1/(n+1)).
$$	

 It is obvious that

$$
 || f(n,\cdot) \ - \ f(\cdot)||G(\psi) \to 0, \ n \to \infty,
$$
but does not exists some function $ \nu(\cdot) $ for which  that $ \nu >> \psi $
and such that

$$
\sup_n ||f(n,\cdot)||G(\nu) < \infty,
$$
as long as $ f(\cdot) \notin G^o(\psi). $ \par

\vspace{3mm}

{\bf Remark 2.}  But the condition {\bf 2} can not be omitted. Let us consider the correspondent example.\par

 Let $ X = [0,1] $ with ordinary Lebesgue measure $ m. $ Let $ \xi = \xi(\omega), \ \omega \in X $ be a standard Gaussian distributed function, \
on the other words,  Normal Standard random variable.\par
 Put $ \psi(p) = \psi_{0.5}(p) =  p^{1/2}, \ p \ge 1. $ It is easy to calculate that $ \xi \in G(\psi_{0.5}). $ Consider the sequence

$$
\xi(n, \omega) = \xi(\omega) \cdot I(|\xi(\omega)| \le n), \ n = 1.2.\ldots.
$$
 We observe:

$$
\xi(n,\cdot) \in GB  \ ( \ = GA = G^o(\psi) \ ),
$$

$$
\forall p \in (1,\infty) \ \Rightarrow \ |\xi(n,\cdot) \ - \ \xi|_p \to 0,
$$
as $ \ n \to \infty; $

$$
 \sup_n ||\xi(n,\cdot)|| G(\psi_{0.5}) < \infty.
$$

But the convergence of $ \xi(n, \cdot) $ to the variable $ \xi(\cdot) $ in the
sense  of $ G(\psi_{0.5}) $ norm is false as long as the limit variable
$ \xi $ does not belongs to the space $ G^o(\psi_{0.5}). $\par

\vspace{3mm}

{\bf Remark 3.} Instead the condition {\bf 1 } it can be presumed that the
compactness of the sequence of the functions $ \{f(n, \cdot)\} $ is true only
for some {\it sequence} of  spaces $ L_{p(i)} $ with
powers $ \{ p(j) \} $ such that

$$
\overline{\lim}_{j \to \infty} p(j) = b; \ \underline{\lim}_{j \to \infty} p(j) = a.
$$
 See in detail \cite{Rao3}, chapter 4,5. \par

\vspace{3mm}

{\bf Remark 4.} The first  condition {\bf 1} in the case of finiteness of measure $ \mu: \ \mu(X) < \infty $ may be replaced on the compactness on the measure, i.e. in the distance

$$
r(f,g):= \inf \{\epsilon, \epsilon > 0, \ \mu \{x: |f(x) \ - \ g(x)|>
\epsilon  \}  < \epsilon \}.
$$

\vspace{4mm}

{\bf B.} In this pilcrow  we consider arbitrary  bounded closed subset $ S $
of the space $ G^o(\psi) = GA(\psi) = GB(\psi): \ S \subset G^o(\psi),
\ \psi \in \Psi(a,b). $ \par

\vspace{3mm}

{\bf Theorem 2.} {\it The closed bounded subset $ S $ of the space $ G^o(\psi) $
or equally in the spaces $ GA(\psi), \ GB(\psi) $
is sequentially compact set if and only if it is compact set in each space
$ L_p, \ p \in (a,b) $ and is equi \ -  \ absolute continuous:} $ S \in EGA. $ \par

\vspace{3mm}

{\bf Proof } is very simple. The sufficient assertion "if"  follows immediately
from the theorem 1 and the known quoted properties of BSGL spaces. The inverse
affirmation may be proved alike the proof of the theorem 2 from the classical
book \cite{Rao1}, chapter 5, section 5.2. \par
 See also \cite{Rao2}, chapters 1,2 and \cite{Krasnoselsky1}, chapters 1,2,3. \par

\vspace{3mm}

 {\bf Remark 5.} The remarks 3 and 4 are true even in the considered case of the sets
 in $ GA(\psi) $ spaces.\par

\vspace{3mm}

{\bf Remark 6.} Suppose the metric space $ (\Sigma, \rho) $ is separable.	
Then for the sets on the space $ G^o(\psi) $ are true the classical criterions
for compactness, belonging to Kolmogorov, Riesz, Schilov,  Frechet etc.\par

\vspace{3mm}

{\bf Remark 6.} Let $ X $ be the closure of open non \ - \ empty bounded set
in the Euclidean  finite \ - \ dimensional space $ R^d $  equipped Lebesgue
measure $  m. $ Then the spaces $ G^o(\psi), \ GA(\psi) $ and $ GB(\psi) $
have a basis, consisting, for instance, from the deformed Haar's functions . \par
 The proof is at the same as in the case of Orlicz spaces, considered, e.g.
in \cite{Krasnoselsky1}, chapter 2, section 12. See also
\cite{Krein1}, chapter 1.\par

\vspace{4mm}

{\bf 4. Applications. Spherical rearrangement of a functions.  } \\

\vspace{3mm}

{\bf 0.} We consider in this section some slight generalization of a main result of Jean van Schaftingen \cite{Schaftingen1} in the $ G(\psi) $ spaces instead
$ L_p $ spaces considered in \cite{Schaftingen1}. \par
 Let $ u: R^d \to R^+ \cup \{ \infty \} $ be a measurable function. A function
$ u^* = u^*(x) = u^*(x; u(\cdot)) $ is called symmetric rearrangement, or
Schwarz symmetrization  of a function $ u, $ iff it is spherical symmetry in
the following sense:  $ \forall \lambda \in R \ \exists r = r(\lambda) \ge 0, $
	
$$
\{x \in R^d: \ u^*(x) > \lambda \} = B(0,r),
$$
where $ B(0,r) $ is the Euclidean ball with center at the origin and
radius $ r; $ \par

$$
\forall \lambda \in R \ m  \{x: u^*(x) > \lambda  \} = m  \{x: u(x) > \lambda  \}.
$$
Recall that $ m $ denotes the Lebesgue measure.\par

This symmetrization there exists for arbitrary function $ u, \  u \in L_p, $  and is unique (up to the set of zero measure).\par
 The symmetrization is used in the study of isoperimetric inequalities, variational problems, theory of partial differential equations etc., see
\cite{Schaftingen1}.\par

\vspace{3mm}

{\bf 1.} As long as

$$
|u^*|_p = |u|_p, \ p \le 1,
$$
we can conclude then  for every function $ \psi \in \Psi(a,b) $

$$
||u^*||G(\psi) = ||u||G(\psi)
$$
(the norm preserving). \\

\vspace{3mm}

{\bf 2.} Since

$$
|u^*(\cdot; u) \ - \ v^*(\cdot; v) |_p  \le |u(\cdot) \ - \ v(\cdot)|_p,
$$
then for every function $ \psi \in \Psi(a,b) $

$$
||u^* \ - \ v^*||G(\psi) \le ||u \ - \ v||G(\psi)
$$
(contraction). \\

\vspace{3mm}

{\bf 3.} By virtue of  the inequality

$$
| \nabla \ u^*|_p \le  | \nabla \ u|_p, \ p \le 1,
$$
we assert then for every function $ \psi \in \Psi(a,b) $

$$
|| \nabla \ u^*||G(\psi) \le  || \nabla  \ u||G(\psi)
$$
(Polya \ - Szeg\"o inequality), if obvious 

$$
\nabla  \ u \in G(\psi).
$$

\vspace{3mm}

{\bf 4.}
 Jean van Schaftingen in \cite{Schaftingen1} offered the consistent iterative
algorithm  ("polarization algorithm") for $ u^*(\cdot) $ computation. He construct the sequence of a functions $ \{ u_n \}, \ u_1 = u,  $  such that
if $ u \in L_p, \ p \ge 1, $ then:

\vspace{3mm}

$$
{\bf A.} | u_n \ - u^* |_p \to 0, \ n \to \infty;
$$

\vspace{3mm}
$$
{\bf B.} |u_n|_p = |u|_p = |u^*|_p;
$$

\vspace{3mm}

{\bf C.} The sequence $ \{ u_n \} $ is compact set in the space $ L_p(R^d). $ \par

 Now we assume that $ u(\cdot) \in G(\psi) $ for some $ \psi \in \ \Psi(a,b),
1 \le a < b \le \infty. $  From the theorems 1  and 2 it \ follows:\par

\vspace{3mm}

{\bf Theorem 3.} {\it For arbitrary function } $ \nu, \nu \in G(\psi), \ 
\nu  >> \psi $

$$
|| u_n \ - \ u^* ||G(\nu) \to 0, \ n \to \infty.
$$

{\it If in addition $ u(\cdot) \in \ G^o(\psi), $ then }

$$
|| u_n \ - \ u^* ||G^o(\psi) = || u_n \ - \ u^* ||G(\psi) \to 0, 
\ n \to \infty.
$$

 Note that by virtue of the property  {\bf 2 } the function $ u^* (\cdot) =
u^*(\cdot,u) $ depended continuously on the source function $ u(\cdot) $
also in the $ G(\psi) $ norm.\par


\end{document}